\documentclass[11pt,letterpaper]{amsart}

\usepackage[letterpaper,margin=1in]{geometry}
\usepackage[all]{xy} 
\usepackage{amsmath, amssymb, amsfonts, latexsym, mdwlist, amsthm, amscd,mathabx,braket}
\usepackage{subfig}
\usepackage{graphicx}
\usepackage{wrapfig}
\usepackage{etex}
\usepackage{tikz-cd}
\usepackage[colorlinks=true, citecolor=blue, urlcolor=blue, linkcolor=blue, pagebackref]{hyperref}

\usepackage{tikz}
\usetikzlibrary{calc,trees,positioning,arrows,chains,shapes.geometric,%
    decorations.pathreplacing,decorations.pathmorphing,shapes,%
    matrix,shapes.symbols}

\tikzset{
>=stealth',
  punktchain/.style={
    rectangle,
    rounded corners,
    draw=black, thick,
    minimum height=3em,
    text centered,
    on chain},
  line/.style={draw, thick, <-},
  eLement/.style={
    tape,
    top color=white,
    bottom color=blue!50!black!60!,
    minimum width=8em,
    draw=blue!40!black!90, very thick,
    text width=10em,
    minimum height=3.5em,
    text centered,
    on chain},
  every join/.style={->, thick,shorten >=1pt},
  decoration={brace},
  tuborg/.style={decorate},
  tubnode/.style={midway, right=2pt},
}

\usepackage{paralist}
\usepackage{enumitem} 
\setdefaultenum{(1)}{(a)}{(i)}{}
\setlist[enumerate,1]{label={\upshape(\arabic*)}}
\setlist[enumerate,2]{label={\upshape(\alph*)},ref=\theenumi\upshape(\alph*)}
\setlist[enumerate,3]{label={\upshape(\roman*)},ref=\theenumi\theenumii\upshape(\roman*)}
\usepackage[capitalize]{cleveref}
\crefname{Prop}{Proposition}{Propositions}
\crefname{Thm}{Theorem}{Theorems}
\crefname{Lem}{Lemma}{Lemmas}
\crefname{enumi}{Case}{Cases}

\def\C{\ensuremath{\mathbb{C}}}

\def\N{\ensuremath{\mathbb{N}}}
\def\P{\ensuremath{\mathbb{P}}}

\def\Z{\ensuremath{\mathbb{Z}}}

\def\codim{\mathop{\mathrm{codim}}\nolimits}

\def\disc{\mathop{\mathrm{disc}}\nolimits}
\def\dim{\mathop{\mathrm{dim}}\nolimits}

\def\im{\mathop{\mathrm{im}}\nolimits}

\def\mod{\mathop{\mathrm{mod}}\nolimits}
\def\min{\mathop{\mathrm{min}}\nolimits}

\DeclareMathOperator{\Pim}{Pim}

\def\rk{\mathop{\mathrm{rk}}}

\def\Span{\mathop{\mathrm{Span}}\nolimits}

\def\o{\ensuremath{\mathfrak{o}}}

\def\MG13{\ensuremath{{\mathcal M}_{\Gamma_1(3)}}}
\def\tildeMG13{\ensuremath{\widetilde{\mathcal M}_{\Gamma_1(3)}}}

\def\into{\ensuremath{\hookrightarrow}}

\newcommand\TFILTB[3]{%
\xymatrix@=1pc{
{0 = {#1}_0} \ar[rr]&&
{{#1}_1} \ar[rr]\ar[ld] &&
{{#1}_2} \ar[r]\ar[ld] &
{\cdots} \ar[r] & { {#1}_{#3-1}} \ar[rr] &&
{{#1}_{#3} = {#1}} \ar[ld]
\\
& *{{#2}_1} \ar@{.>}[ul] &&
{{#2}_2} \ar@{.>}[ul] & &&&
{{#2}_{{#3}}} \ar@{.>}[ul]
}}

\makeatletter
\newtheorem*{rep@theorem}{\rep@title}
\newcommand{\newreptheorem}[2]{%
\newenvironment{rep#1}[1]{%
 \def\rep@title{#2 \ref{##1}}%
 \begin{rep@theorem}}%
 {\end{rep@theorem}}}
\makeatother

\newtheorem{Thm}{Theorem}[section]
\newreptheorem{theorem}{Theorem}
\newtheorem{proposition}[Thm]{Proposition}

\newtheorem{lemma}[Thm]{Lemma}

\newtheorem{corollary}[Thm]{Corollary}
\newreptheorem{Cor}{Corollary}

\newreptheorem{Con}{Conjecture}

\newtheorem*{theorem*}{Theorem}
\newtheorem*{lemma*}{Lemma}
\newtheorem*{proposition*}{Proposition}
\newtheorem*{conjecture*}{Conjecture}
\newtheorem*{corollary*}{Corollary}
\newtheorem*{problem*}{Problem}

\newtheorem{Thm-int}{Theorem}

\theoremstyle{definition}
\newtheorem{Def-s}[Thm]{Definition}
\newtheorem{definition}[Thm]{Definition}

\newtheorem{example}[Thm]{Example}

\def\C{\ensuremath{\mathbb{C}}}

\def\N{\ensuremath{\mathbb{N}}}
\def\P{\ensuremath{\mathbb{P}}}

\def\Z{\ensuremath{\mathbb{Z}}}

\def\CC{\ensuremath{\mathcal C}}

\def\HH{\ensuremath{\mathcal H}}
\def\II{\ensuremath{\mathcal I}}
\def\JJ{\ensuremath{\mathcal J}}

\def\SS{\ensuremath{\mathcal S}}

\def\ZZ{\ensuremath{\mathcal Z}}

\newcommand{\ignore}[1]{}
\begin{document}

\title{Supporting rank and the intersection of all Hassett divisors}

\author{Elad Gal}
\address{Faculty of Mathematics \\
Technion, Israel Institute of Technology}
\email{galelad@campus.technion.ac.il}
\author{Howard Nuer}
\address{Faculty of Mathematics \\
Technion, Israel Institute of Technology}
\email{hnuer@technion.ac.il}

\thanks{}

\begin{abstract}
We prove that the dimension of the intersection $\ZZ$ of all Hassett divisors of special cubic fourfolds is sixteen.  We do this by studying which subsets of the natural numbers $\N$ can be obtained as the image of a positive-definite integral quadratic form and what the minimal possible rank of such a form is.  In particular, for the subset of $\N$ consisting of all possible discriminants of special cubic fourfolds, we show this rank is four and that this is the codimension of $\ZZ$ in $\CC$, the twenty-dimensional moduli space of cubic fourfolds.    
\end{abstract}

\maketitle
\setcounter{tocdepth}{1}
\tableofcontents
\section{Introduction}\label{Sec:Intro}
Cubic fourfolds occupy an important place in modern algebraic geometry for a number of reasons: the elusiveness and difficulty of their rationality problem, their relation to K3 surfaces, and the relative manageability of their Hodge theory, to name a few. 
This latter point makes studying their moduli very approachable using tools from lattice theory. 
Indeed, for a cubic fourfold $X$, a smooth degree three hypersurface of $\P^5_\C$, its algebraic cohomology $A(X):=H^4(X,\Z)\cap H^{2,2}(X)$ is a positive-definite lattice satisfying useful properties.  
It always contains the square of the hyperplane class, $h_X^2$, which is distinguished in the sense of Yang and Yu \cite{YangYu}, that is, it has norm $3$ and even orthogonal complement (see \cref{def:distinguishedelement}).  
For a very general cubic fourfold, we have $A(X)=\Z h_X^2$, and these cubic fourfolds are the first among those conjectured to be irrational.  

By building on Voisin's Torelli theorem for cubic fourfolds \cite{Voisin:Cubics},  Hassett pioneered the study of special cubic fourfolds, those with $\rk A(X)\geq 2$, in his seminal work \cite{Has00}.  Hassett proved that in the twenty-dimensional moduli space $\CC$ of cubic fourfolds the locus of special cubic fourfolds is the infinite union of irreducible divisors $\CC_d$, which denotes the moduli space of special cubic fourfolds of discriminant $d$.  Moreover, he showed that the Hassett divisor $\CC_d$ is nonempty if and only if 
$$d\in\HH:=\{d\in\N | d>6, d\equiv 0,2\pmod 6\},$$
which we call the Hassett subset.  

Generalizing the ideas of Voisin and Hassett even further, Yang and Yu introduced the notions of an admissible positive-definite lattice $M$ with a distinguished element and of an $M$-polarizable cubic fourfold \cite{YangYu}.  They also gave precise conditions for the moduli space $\CC_M$ of $M$-polarizable cubic fourfolds to be irreducible and non-empty (see \cref{theorem:YangYu} below).  Historically, some of the first examples in this context of admissible lattices of rank larger than two were obtained by considering $A(X)$ for $[X]$ in the intersection of two or more Hassett divisors.  This, among other reasons, led Yang and Yu to study the intersection of all Hassett divisors,
$$\ZZ:=\bigcap_{d\in\HH}\CC_d.$$
They showed that the Fermat cubic fourfold lies in $\ZZ$ so that it is nonempty \cite[Theorem 1.2]{YangYu} and proved the estimate \cite[Remark 7.15]{YangYu}
$$13\leq\dim\ZZ\le16,$$
leaving the determination of the dimension of $\ZZ$ as an open question \cite[Question 7.14]{YangYu}.

The main result of this note answers that open question:
\begin{Thm}[{See \cref{Thm:MainTheorem}}]
    The intersection $\ZZ$ of all Hassett divisors has dimension $16$.
\end{Thm}
With the powerful tool of Voisin's Torelli Theorem in hand, we prove this result by translating it into a question in lattice theory that generalizes some of the most classical and beautiful results in the arithmetic side of the subject.  
Indeed, we were lead to consider the following questions in the arithmetic theory of quadratic forms: which subsets $S\subset\N$ can be obtained as the image of a positive-definite integral quadratic form, and for such a subset $S$ what is the smallest rank of a positive-definite integral quadratic form with image $S$? 
We call such subsets \emph{formally-supported} and denote the minimal rank of a supporting quadratic form by $r(S)$, which we call the \emph{supporting rank} of $S$.
For example, a classical notion in the arithmetic theory of integral quadratic forms is that of a universal positive-definite quadratic form, those whose image is $\N$.  
As such forms exist, $\N$ is formally-supported, and it follows from Lagrange's four squares theorem and \cref{prop:ternaryforms} below that the supporting rank of $\N$ is four. 
Furthermore, we can ask if any element in a subset $S$ formally-supported by a positive-definite quadratic form is the image of a primitive vector.  We say such subsets $S\subset\N$ are \emph{primitively-formally-supported} and call the minimal rank of a primitively supporting quadratic form the \emph{primitive supporting rank} of the subset, denoted $r_{prim}(S)$.
Stated in these terms, our main result is equivalent to the following:
\begin{Thm}[{See \cref{corollary:RankHFour}}]
    The Hassett subset $\HH$ satisfies $r_{prim}(\HH)=4$.
\end{Thm}

While we have determined in this note the maximal dimension of an irreducible component of $\ZZ$, it remains unclear if $\ZZ$ is equidimensional.  
That is, are there higher codimensional irreducible components?
A related but independent question is: how many irreducible components does $\ZZ$ have?
In the course of our work here, we have found at least two irreducible components of $\ZZ$ of dimension sixteen. 
In the lattice theoretic language above, we can ask more generally if a lattice supporting a given subset $S\subset\N$ contains a sublattice of the minimal rank $r(S)$ that still supports $S$.
Furthermore, we can ask how many distinct non-isometric lattices that support $S$ are minimal (under inclusion) with respect to this support property. 
We leave these questions to subsequent work.

\section{Lattices, formally-supported subsets, and supporting rank}
\subsection{Lattices and integral quadratic forms}
A \emph{lattice} is a finite rank free $\Z$-module $M$ together with a symmetric and non-degenerate bilinear form $\langle-,-\rangle_M\colon M\times M\to\Z$.  
The \emph{discriminant} of a lattice $M$, denoted $\disc(M)$, is the determinant of the Gram matrix with respect to an arbitrary basis of $M$.  
An element $v\in M$ is called a \emph{root} if $\langle v,v\rangle_M=2$.
A lattice $M$ is called \emph{even} if $\langle v,v\rangle_M\in 2\Z$ for all $v\in M$, and \emph{odd} otherwise.  

A sublattice $N\subset M$ is called \emph{primitive} if $M/N$ is torsion-free.  
We denote by $N^\perp$ the \emph{orthogonal complement} of $N$, that is, the set $\{v\in M | \langle v,N\rangle_M=0\}$.

An \emph{integral quadratic form} is a function $q\colon M\to\Z$ on a finite rank free $\Z$-module $M$ such that
\begin{enumerate}
    \item $q(n v)=n^2q(v)$ for all $n\in\Z$ and $v\in M$, and 
    \item the function $Q\colon M\times M\to\Z$ given by $Q(u,v)=q(u+v)-q(u)-q(v)$ is bilinear and non-degenerate.  
\end{enumerate}
An integral quadratic form $q$ is said to be \emph{integer-matrix} if the associated function $Q$ satisfies $2\mid Q(u,v)$ for all $u,v\in M$.  This is equivalent to $\frac{1}{2}Q\colon M\times M\to\Z$ defining the structure of a lattice on $M$ such that $q(v)=\frac{1}{2}Q(v,v)$ for all $v\in M$.

We say a lattice $(M,\langle-,-\rangle_M)$ (resp. integral quadratic form $(M,q)$) is positive (or positive-definite) if $\|v\|_M:=\langle v,v\rangle_M>0$ (resp. $q(v)>0$) for all $0\ne v\in M$.  
We denote the image $q(M)\subset\N$ of a positive integral quadratic form by $\im(q)$ or $\im(M)$ (by abuse of notation).
\subsection{Formally-supported subsets}
The celebrated $15$ and $290$ Theorems give precise, easily verified criteria for a positive integer-matrix (resp. integral) quadratic form to be universal; that is, to have image equal the entire set of natural numbers.  A number of natural questions emerge from these types of results.  For example, one can ask for what subsets \(S\subset\N\) does there exist a positive integral (or integer-matrix) quadratic form $f$ such that $\im(f)=S$?  We can also ask further for the subset $S$ to be represented by primitive vectors.  These considerations lead to the following definitions:
\begin{definition}
A subset $S\subset\N$ is called \textbf{formally-supported} if there exists a positive integral quadratic form $f$ with $\im(f) = S$.  Similarly we say $S$ is \textbf{integrally-formally-supported} if we we may take $f$ to be  integer-matrix, and we say that $S$ is \textbf{primitively-formally-supported} if $S$ is formally-supported (say by $f$) in such a way that for every $n\in S$, $n=f(v)$ has a primitive solution $v\in\Z^n$.
\end{definition}

For example, Lagrange's four squares theorem shows that \N\ is integrally-formally-supported. The $15$-Theorem \cite{Conway,Bhargava} (respectively, the $290$-Theorem \cite{Hahn:290}) can be rephrased to say that $\N\setminus\{n\}$ is integrally-formally-supported (respectively, formally supported) if and only if $s$ is one of finitely many values, the largest of which is $15$ (respectively, $290$). 
Going one step further, \cite[Theorem 1]{AlmostUniversal} says that $\N\setminus\{n,m\}$ is integrally-formally-supported if and only if $\{n,m\}$ is one of $73$ pairs of natural numbers with $n<m$. It is worth noting that integrally-formally-supported implies formally-supported, but not vice-versa.  Indeed, $\N\setminus\{290\}$ is formally-supported but not integrally-formally-supported.

\subsection{Supporting rank}
When studying a formally-supported set it is also natural to ask how efficiently its elements can be represented by a positive integral quadratic form.
So we define the rank of a formally-supported set $S\subset\N$ as follows.
\begin{definition}
The \textbf{supporting rank} of $S$ is defined to be
$$r(S):=\min\left\{\rk(f)|\im(f)=S\right\},$$
where we denote by $\rk(f)$ the number of variables appearing in $f$.
Similarly, we can define $r_{int}$ if we consider only integer-matrix positive forms and $r_{prim}$ if we require the solution to be a primitive element.
\end{definition}

Our main theorem, when described in the language of quadratic forms, is a statement about the supporting rank and primitive supporting rank of a certain set of natural numbers that arises in an algebraic geometry problem. 
We begin here with a classical result that ternary forms necessarily miss many numbers.
\begin{proposition}\label{prop:ternaryforms}
For any positive ternary rational quadratic form $f$, there is an infinite set \(S=\{s_i\}_{i \in \mathbb{N}}\) of square-free integers which is disjoint from \(f(\mathbb{Q}^{\oplus3})\)
\end{proposition}
\begin{proof}
The above result also follows from \cite[p. 142]{Conway:SensualQuadratic} or \cite[Theorem 1.3]{TernaryForms}.
\end{proof}

This theorem shows that any form of rank $\leq3$ misses infinitely many natural numbers and implies the following result:
\begin{corollary}\label{cor:RankAtLeastOne}
If $S\subset\N$ is a subset with $r(S)\leq3$, then $|\N\setminus S|=\infty$. This is true with regard to $r_{int}$ and $r_{prim}$ as well.
\end{corollary}

\section{The Hassett subset}
The main subset of \N\ of interest to us here is the \textbf{Hasset subset} $\HH$ defined by
$$\HH:=\{d\geq8|d\equiv0,2\pmod 6\},$$
which is the set of $d$ for which the Hassett divisor $\CC_d$ of special cubic fourfolds of discriminant $d$ is nonempty, as mentioned in \cref{Sec:Intro}.

We study the question of representing this set as the image of a positive integral quadratic form. 
We note first that $\HH$ is easily seen to be primitively-integrally-formally-supported as the image of the integer-matrix positive quadratic form
$$f=8x^2+12y^2+6xy+18z^2+12(a^2+b^2+c^2+d^2).$$

In fact, we show in \cref{Sec:SupportingRankH} that $\HH$ is primitively-integrally-supported by a positive integer-matrix quadratic form of rank $4$.
For this, we begin by recalling an important definition from lattice theory:

\begin{definition}\label{def:distinguishedelement}
If $M$ is a lattice with positive pairing $\langle-,-\rangle_M$, an element $\o$ will be called a \textbf{distinguished element} if $\langle\o,\o\rangle_M=3$ and $\o^\perp\subset M$ is an even sublattice.
\end{definition}
\begin{definition}
For any rank $n$ lattice $M$ with positive pairing and distinguished element $\o$,  we will denote by $D(M)$ the rank $n-1$ lattice $M/\Z\o$ with associated quadratic form defined by 
$$\|v+\Z\o\|_{D(M)}=\disc(\Span(\o,v))_M.$$
\end{definition}
One can check that if $M$ has a distinguished element $\o$ and contains no roots (that is, elements of square two), then $\im(D(M))\subset\HH$.
Surprisingly, the converse also holds:

\begin{proposition}
For any rank $n$ positive-definite integer-matrix form $f$ with $f(\Z^{\oplus n})\subset\HH$ there exists a rank $(n+1)$ positive-definite lattice $M$ with a distinguished element $\o$ and containing no roots such that
$$ f(v) = \disc(\Span(\o,v))_M = 3\|v\|_M-\langle\o,v\rangle_M^2 $$
\end{proposition}
\begin{proof}
Take an integral basis $\{v_1,v_2,...,v_n\}$ for the rank $n$ lattice $L$ with pairing $\|v\|_L=f(v)$. We know that for any $v\in L$ we have $||v||_L\equiv0,2\pmod6$. If $\|v_{n-1}\|_L\equiv\|v_n\|_L\equiv2\pmod 6$, then 
$$\|v_{n-1}\pm v_n\|_L = \|v_{n-1}\|_L+\|v_n\|_L\pm2\langle v_{n-1},v_n\rangle_L \equiv 4\pm2\langle v_{n-1},v_n\rangle_L\pmod 6.$$
This norm is in $\HH$, so $\langle v_{n-1},v_n\rangle_L\not\equiv0(3)$, and thus one of the vectors $v_{n-1}\pm v_n$ has norm divisible by 6. Replacing $v_n$ with this vector, we get a new basis with $||v_n||_L\equiv0\pmod 6$.  Repeating this inductively, we can assume our basis $\{v_1,\dots,v_n\}$ satisfies $||v_i||_L\equiv0\pmod 6$ for any $2\leq i\leq n$.

Let us consider $||v_i\pm v_j||_L\in\HH$.  If both have norm divisible by $6$ then 
$$\|v_i\pm v_j\|_L=\|v_i\|_L+\|v_j\|_L\pm2\langle v_i,v_j\rangle_L\equiv\pm2\langle v_i,v_j\rangle_L\pmod 6.$$
If $\langle v_i,v_j\rangle_L\equiv 1\pmod 3$, then $\|v_i- v_j\|_L\equiv 4\pmod 6$ which is not in $\HH$.  Similarly, if $\langle v_i,v_j\rangle_L\equiv 2\pmod 3$, then $\|v_i+ v_j\|_L\equiv 4\pmod 6$ which is not in $\HH$.  Thus we can conclude that  $\langle v_i,v_j\rangle_L\equiv0\pmod 3$.  The same reasoning shows that this holds if one of them has norm equivalent to $2\pmod 6$.

Now, define the new rank $(n+1)$ positive-definite lattice $M$ with distinguished element $\o$ to be the lattice with basis $\{\o,v_1',\dots,v_n'\}$ and pairing:

$$\begin{pmatrix}
M  &\o & v_1' &\dots& v_i'&\dots\\
\o & 3 &  i  &\dots&0&\dots \\
v_1'& i &  \frac{||v_1||_L+i}{3}   &\dots&  \frac{\langle v_i,v_1\rangle_L}{3} &\dots \\
\vdots&\vdots&\vdots&\vdots&\vdots&\vdots\\
v_j'&0&  \frac{\langle v_1,v_j\rangle_L}{3}&\dots  & \frac{\langle v_i,v_j\rangle_L}{3}&\dots  \\
\vdots&\vdots&\vdots&\vdots&\vdots&\vdots\\
\end{pmatrix},$$
where $i=\begin{cases}
    0&\text{if } f(v_1)\equiv0\pmod 6 \\
    1&\text{if } f(v_1)\equiv2\pmod 6\\
\end{cases}$.  
Writing $v=\sum_{i=1}^n c_iv_i'$, we get
\begin{align*}
  &  3\|v\|_M-\langle\o,v\rangle_M^2=3\sum_{i=1}^n\sum_{j=1}^n\langle c_iv_i',c_jv_j'\rangle_M-\sum_{i=1}^n\sum_{j=1}^n\langle\o,c_iv_i'\rangle_M\langle\o,c_jv_j'\rangle_M\\
&=3\left(\frac{ic_1^2}{3}+\sum_{i=1}^n\sum_{j=1}^n\frac{\langle c_iv_i,c_jv_j\rangle_L}{3}\right)-c_1^2i^2
=ic_1^2+\sum_{i=1}^n\sum_{j=1}^n\langle c_iv_i,c_jv_j\rangle_L-c_1^2i^2=\\
&=i(1-i)c_1^2+\langle\sum_{i=1}^nc_iv_i,\sum_{j=1}^nc_jv_j\rangle_L=i(1-i)c_1^2+\|\sum_{i=1}^n c_iv_i\|_L=f(v)+i(1-i)c_1^2
\end{align*}
But $i\in\{0,1\}$, so $i(1-i)=0$. 
Thus $f(v)=3\|v\|_M-\langle\o,v\rangle_M^2$, as required.
\end{proof}
\begin{corollary}\label{Cor:RankAtLeast4}
$r(\HH)\geq4$.
\end{corollary}
\begin{proof}
Let $L$ be a rank $n$ lattice with $\im(L)=\HH$.  From the proposition, there exists a  rank $(n+1)$ lattice $M$ with $D(M)=L$. This implies that for any $d\geq2$ there exist a vector $v\in M$ such that $6d=\disc(\Span(\o,v))$. Now $\Span(\o,v)$ is a rank $2$ sublattice with Gram matrix
$$\begin{pmatrix}
- & \o & v\\
\o&  3 &  \langle\o,v\rangle_M \\
v& \langle\o,v\rangle_M & \|v\|_M \\
\end{pmatrix},$$
so 
$$6d=3\|v\|_M-\langle\o,v\rangle_M^2$$
implies that $\langle\o,v\rangle_M$ is dividible by 3.
Replacing $v$ by $v'=v-\frac{\langle\o,v\rangle_M}{3}\o$, we get a basis $\{\o,v'\}$ for this sublattice such that $\langle\o,v'\rangle_M=0$. Then $6d=\disc(\Span(\o,v'))$ means that $\|v'\|_M=2d$.

So, for any $d\in\N\setminus\{1\}$, there exist $v\in\o^\perp\subset M$ with norm $2d$. But $$\o^\perp\subset M$$ is even, so dividing the intersection form on $\o^\perp$ by $2$ we get a lattice with vectors of any norm except $1$. Thus this $n$-form supports $\N\setminus\{1\}$. It follows from \cref{cor:RankAtLeastOne} that $n\geq4$.
\end{proof}

\section{The supporting rank of $\HH$}\label{Sec:SupportingRankH}
We showed in \cref{Cor:RankAtLeast4} that $r_{prim}(\HH),r(\HH)\geq 4$.  In this section we will prove that 
$$r_{prim}(\HH)=r(\HH)=4.$$ To do this it suffices to find one lattice of rank $4$ whose primitive image is $\HH$.

Let $M$ be the following rank $5$ lattice containing a distinguished element and having Gram matrix
$$\begin{pmatrix}
M & \o & a & b & c & d \\
\o & 3 & 1 & 0 & 0 & 0 \\
a & 1 & 3 & 2 & 0 & 2 \\
b & 0 & 2 & 6 & 2 & 0 \\
c & 0 & 0 & 2 & 4 & 0 \\
d & 0 & 2 & 0 & 0 & 4 \\
\end{pmatrix}.$$
As this lattice has a distinguish element, we may apply the definition of $D(M)$ to get a rank $4$ lattice with Gram matrix
$$D(M) = \begin{pmatrix}
D(M) & a & b & c & d \\
a & 8 & 6 & 0 & 6 \\
b & 6 & 18 & 6 & 0 \\
c & 0 & 6 & 12 & 0 \\
d & 6 & 0 & 0 & 12 \\
\end{pmatrix}$$
and corresponding rank $4$ integer-matrix quadratic form
\begin{equation}\label{eq:originalform}
    f(x,y,z,u)=8x^2+12xy+18y^2+12yz+12z^2+12u^2+12xu.
\end{equation}

\begin{Thm}\label{Prop:MainComponentDimension}
Let $f$ be the integer-matrix quadratic form defined in \eqref{eq:originalform}.  Then
\[\Pim(f)=\im(f)=\HH\]
\end{Thm}
\begin{proof}
We note first that we can write \eqref{eq:originalform} as
\begin{equation}\label{eq:positiveform}
     f(x,y,z,u)=\frac{13}{5}x^2+15\left(y+\frac{2}{5}x\right)^2+3\left(y+2z\right)^2+3\left(x+2u\right)^2  
\end{equation}
so that \(f\) is positive-definite.  Moreover, $f(x,y,z,u)\equiv 2x^2\pmod 6$, so any element in $\im(f)$ must be $0,2\pmod 6$.
Now we must eliminate $2$ and $6$ from the image.  Observe first that if $2=f(x,y,z,u)$, then from \eqref{eq:positiveform} we would have $2\geq\frac{13}{5}x^2$ which forces $x=0$, and from \eqref{eq:originalform} it would follow that $6|f(0,y,z,u)=2$, a contradiction.  Similarly, if $6=f(x,y,z,u)$, then since every term in \eqref{eq:originalform} but $8x^2$ is divisible by $3$, it would follow that $3|x$ so that from \eqref{eq:positiveform} we would have 
$$6\geq\frac{13\cdot 9}{5}\left(\frac{x}{3}\right)^2,$$ forcing $x=0$.  It then follows further from \eqref{eq:positiveform} that $y$ and $u$ would have to be $0$ and then that $z=0$, a contradiction.  
This proves that $\Pim(f)\subset\im(f)\subset\HH$.

To finish the proof it suffices to prove that $\HH\subset\Pim(f)$.  To that end, we write
$$5f(x,y,z,u)=13x^2+3(5y+2x)^2+15(y+2z)^2+15(x+2u)^2$$

We split the values of $d\in\HH$ into cases  based on their residue modulo $24$ and the shape of the primitive solution we find in each case as follows, where primitivity of the solution is automatic or easily shown by the construction:
\begin{itemize}
    \item $\{8,14,20\}+24\N_0\subset f(\{1\}\times\Z\times\Z\times\Z)$
    \item $\{26\}+24\N_0\subset f(\{2\}\times\Z\times\Z\times\Z)$
    \item $\{12,24,30\}+24\N_0\subset f(\{3\}\times((\Z\times\Z\times\Z)\setminus(3\Z\times3\Z\times3\Z)))$
    \item $\{18\}+24\N_0\subset f(\{6\}\times((\Z\times\Z\times\Z)\setminus(3\Z\times3\Z\times3\Z)))$
\end{itemize}

\noindent\textbf{Case 1: $\{8,14,20\}+24\N_0$}

We will show that$$24\N_0+\{8,14,20\}\subset f(\{1\}\times\Z\times\Z\times\Z)$$
    As the first coordinate is $1$, it is automatic in this case that a solution is primitive.
    Let $n=24d+14$ or $n=12d+8$ be an element of this set. We need to find integers $y,z,u\in\Z$ such that:
    $$120d+70=5n=5f(1,y,z,u)=13+3(5y+2)^2+15(y+2z)^2+15(1+2u)^2,\quad\text{or}$$
$$60d+40=5n=5f(1,y,z,u)=13+3(5y+2)^2+15(y+2z)^2+15(1+2u)^2,$$
respectively.  This is equivalent to solving in integers the equations
    $$40d+19=(5y+2)^2+5(y+2z)^2+5(1+2u)^2,\quad\text{or}$$
$$20d+9=(5y+2)^2+5(y+2z)^2+5(1+2u)^2,$$
respectively.  

In case $n=24d+14$, we can see easily that $40d+19\equiv3\pmod 8$. So, there exist $a,b,c\in\Z$ odds such that $$40d+19=a^2+b^2+c^2$$
    We know that $a^2,b^2,c^2\pmod 5\in\{0,\pm1\}$. As $40d+19\equiv4\pmod 5$, they cannot all be equivalent to $0$ or 1 modulo $5$. Thus we can assume without loss of generality that $a^2\equiv-1\equiv4\pmod 5$. Reducing modulo $10$ we see that up to replacing $a$ by $-a$, we can write $a=5(2a'+1)+2$ for some $a'\in \Z$. 
Similarly, in case $n=12d+8$, we see that $20d+9\equiv1\pmod 4$, so there exist three integers $a,b,c$, two even and the third odd, such that 
    $$20d+9=a^2+b^2+c^2.$$
    As in the previous case, we can assume without loss of generality that $a\equiv2\pmod 5$. Now, $20d+9-a^2=b^2+c^2$ is sum of two squares. By \cref{proposition:Lowering the 5-ht} we can conclude that we can write this value as a sum of two squares $b^2+c^2$ such that $b\equiv2(5),\ c\equiv1(5)$.
    Now, we have three integers numbers $a,b,c$ such that $a^2+b^2+c^2=20d+9$, such that two of them are equivalent to $2$ modulo $5$,, and such that two of them are even. Using the pigeonhole principle we can conclude that at least one of the numbers must satisfy both of these latter conditions. Thus without loss of generality, we can assume that $a=10a'+2$ for some $a'\in\Z$.
    
    In summation, in either case we have integers $a',b,c$ such that 
    $$40d+19=(5(2a'+1)+2)^2+b^2+c^2,\quad\text{or}$$
    $$20d+9=(10a'+2)^2+b^2+c^2,$$
    respectively.
    As $40d+19-(5(2a'+1)+2)^2$, respectively $20d+9-(10a'+2)^2$, is divisible by $5$, it follows from \cref{lemma:Lagrange5} that there exists integers $b',c'$ such that 
    $$40d+19-(5(2a'+1)+2)^2=b^2+c^2=5b'^2+5c'^2,\quad\text{or}$$
     $$20d+9=(10a'+2)^2+5b'^2+5c'^2,$$
     respectively.  
     In the first case, by computation modulo $4$ we see that $b',c'$ are odd, so we can conclude that there exist integers $a',b'',c''$ such that
    $$40d+19=(5(2a'+1)+2)^2+5(1+2b'')^2+5(1+2c'')^2$$
    Using $\bar{b}:= b''-a'$ we can rewrite 
    $$40d+19=(5(2a'+1)+2)^2+5(2a'+1+2\bar{b})^2+5(1+2c'')^2$$
    Finally, substituting $(y,z,u):=(2a'+1,\bar{b},c'')$ we find the desired solution.
    In the second case, we can assume without loss of generality that $b'$ is even and $c'$ is odd. Defining $b'':=\frac{b'-2a'}{2},\ c'':=\frac{c'-1}{2}$ we get 
    $$20d+9=(10a'+2)^2+5(2a'+2b'')^2+5(1+2c'')^2$$   
    Finally, substituting $(y,z,u)=(2a',b'',c'')$ we find the desired solution.

\noindent\textbf{Case 2: $\{26\}+24\N_0$} 

   We will show that
   $$24\N_0+\{26\}\subset f(\{2\}\times\Z\times\Z\times\Z)$$
     Let $n=24d+26$ be an element of this set. We need to find integers
     $y,z,u\in\Z$ such that:
     $$120d+130=5n=5f(2,y,z,u)=52+3(5y+4)^2+15(y+2z)^2+15(2+2u)^2$$
     Observe immediately that any solution must be primitive since it could only be divisible by $2$ but $130$ is not divisible by $4$.  Solving this equation is equivalent to solving in integers the equation:
     $$40d+26=(5y+4)^2+5(y+2z)^2+5(2+2u)^2$$
     As $40d+26\equiv 2\pmod 8$, it can be represented as a sum of three squares $a^2+b^2+c^2$, one even and two odd. By computation modulo $5$, we see that without loss of generality we may assume that $a^2\equiv 1\pmod 5$ so that up replacing $a$ by $-a$, we may assume $a\equiv -1\equiv4\pmod 5$.  Thus
     $$40d+26=(5a'+4)^2+b^2+c^2$$
     for some $a'\in\Z$. Using \cref{proposition:Lowering the 5-ht}, we can rewrite it as 
     $$40d+26=(5a'+4)^2+(5b'+4)^2+(5c'+2)^2$$
     for some $b',c'\in\Z$.
     We know that two of the three numbers are odd, so we can assume without loss of generality that 
     $5a'+4=10a''+9$ for some $a''\in\Z$.
     Now, applying \cref{lemma:Lagrange5} to $(5b'+4)^2+(5c'+2)^2$ we find that there exist $a'',b'',c''\in\Z$ such that
     $$40d+26=(10a''+9)^2+5(1+2b'')^2+5(2+2c'')^2$$
     Finally, substituting $(y,z,u)=(2a''+1,b''-a'',c'')$ we find the desired solution.

\noindent\textbf{Case 3: $\{12,24,30\}+24\N_0$}  

As $f(0,0,0,1)=12$, it suffices to show that
$$24\N_0+\{24,30,36\}\subset f(\{3\}\times\Z\times\Z\times\Z)$$
     with primitive input.
     Let $n=24d+24,30,36$ be an element of this set. We need to find integers $y,z,u\in\Z$, \textbf{not all divisible by $3$}, such that
$$120d+\begin{cases}
    120\\
    150\\
    180\\
\end{cases}=5n=5f(3,y,z,u)=13\cdot9+3(5y+6)^2+15(y+2z)^2+15(3+2u)^2$$
     This is equivalent to solving in integers (not all divisible by 3) the equation:
     $$40d+\begin{cases}
    1\\
    11\\
    21\\
\end{cases}=(5y+6)^2+5(y+2z)^2+5(3+2u)^2$$

Similarly to the first cases, $LHS$ can be represented as a sum of threes squares $a^2+b^2+c^2$, one of them odd, and the others equivalent to $\frac{LHS-1}{2}$ modulo $2$. 
By computation modulo $5$, we may assume that at least one of them is equivalent $1\pmod 5$, say $a\equiv1\pmod 5$. Thus we have integers $a',b,c\in\Z$ such that
$$40d+\begin{cases}1\\11\\21\end{cases}=(5a'+1)^2+b^2+c^2$$
Using \cref{proposition:Lowering the 5-ht} we can rewrite it as 
$$40d+\begin{cases}1\\11\\21\end{cases}=(5a'+1)^2+(5b'+1)^2+(5c'+3)^2$$
for some $b',c'\in\Z$.
We know that one of the three squared numbers is odd and the other two have the same parity as $\frac{LHS-1}{2}$, so without loss of generality we may assume that $a\equiv 1\pmod 5$ and $a\equiv \frac{LHS-1}{2}\pmod 2$.
This implies that there exists an integer $a''$ such that
     $$40d+\begin{cases}1\\11\\21\end{cases}=(5(2a''+\frac{LHS-1}{2})+6)^2+(5b'+1)^2+(5c'+3)^2$$
     Using \cref{lemma:Lagrange5}, we can find $b'',c''\in\Z$ such that 
     $$40d+\begin{cases}1\\11\\21\end{cases}=(5(2a''+\frac{LHS-1}{2})+6)^2+5b''^2+5c''^2.$$
     Reducing modulo $4$ shows that one of $b'',c''$ must be odd while the other has the same parity as $\frac{LHS-1}{2}$.  Let's assume that $b''\equiv\frac{LHS-1}{2}\pmod$, so that $c''$ is certainly odd.  Then substituting $(y,z,u)=(2a''+\frac{LHS-1}{2},\frac{b''-\frac{LHS-1}{2}}{2}-a'',\frac{c''-3}{2})$ we find the desired solution. 
     The existence of a primitive solution now follows from \cref{proposition:Mod3Solutions} since $LHS-5c''^2$ is even and equal to $(5(2a''+\frac{LHS-1}{2})+6)^2+5b''^2$ so we may replace a solution divisible by $3$ by one with $y,z$ not divisible by $3$.

\noindent\textbf{Case 4: $\{18\}+24\N_0$}
We will show that
$$24\Z+\{18\}\subset f(\{6\}\times\Z\times\Z\times\Z)$$
     primitively. 
     Let $n=24d+18$ be an element of this set. We need to find integers $y,z,u\in\Z$, \textbf{not all divisible by $3$} such that:
     $$120d+570=5n=5f(6,y,z,u)=13\cdot36+3(5y+12)^2+15(y+2z)^2+15(6+2u)^2$$
     This is equivalent to solving in integers (not all divisible by $3$) the equation:
     $$40d+34=(5y+12)^2+5(y+2z)^2+5(6+2u)^2.$$  
As $40d+34\equiv2\pmod 8$, it can be written as a sum of three squares 
    $$40d+34=a^2+b^2+c^2,$$
and by computation modulo $4$ we see that exactly two of them are odd and the last one is even. As $40d+34\equiv -1\pmod 5$ and $a^2,b^2,c^2\equiv0,\pm1\pmod 5$, it follows that one of the squares must be $-1$, so without loss of generality we may assume that $a\equiv2\pmod 5$. 

As usual we consider $40d+34-a^2$ which is divisible by $5$ and equal to $b^2+c^2$.  Applying \cref{proposition:Lowering the 5-ht} we can choose $b$ and $c$ so that $b^2\equiv-c^2\equiv1\pmod 5$.
Up to replacing $b$ by $-b$ and $c$ by $-c$ we can assume $b\equiv 1\pmod 5$ and  $c\equiv 2\pmod 5$.  
But now $a\equiv c\equiv 2\pmod 5$, and exactly two out of $\{a,b,c\}$ are odd. Thus up to switching $a$ and $c$, we can assume that $a\equiv 7\pmod 10$, so there exists $a'\in\Z$ such that 
$$40d+34=(10a'+7)^2+b^2+c^2.$$
Using \cref{lemma:Lagrange5} and the fact that $40d+34-(10a'+7)^2\equiv0\pmod 5$ we get that there exist integers $b',c'$ such that
$$40d+34=(10a'+7)^2+5b'^2+5c'^2$$
By computation modulo $2$ we see that there exist integers $b'',c''$ such that
$$40d+34=(10a'+7)^2+5(2b''+1)^2+5(2c'')^2.$$
Then substituting $(y,z,u)=(2a'-1,b''-a'+1,c''-3)$ we find the desired solution. 
The existence of a primitive solution now follows from \cref{proposition:Mod3Solutions} since $40d+34-5(2c'')^2$ is even and equal to $(10a'+7)^2+5(2b''+1)^2$ so we may replace a solution divisible by $3$ by one with $y,z$ not divisible by $3$.  

As we have shown that any element of $\HH\ $ can be represented by a primitive vector in $\Z^4$, it follows that 
$$\Pim(f)=\im(f)=\HH.$$
\end{proof}
This gives us the important result
\begin{corollary}\label{corollary:RankHFour}
$r(\HH)=r_{prim}(\HH)=4$.
\end{corollary}
\section{Lattice polarizable cubic fourfolds and their moduli}
Now we turn toward our goal of applying this lattice theoretic result to the problem in algebraic geometry that inspired our work here.
\subsection{Definitions}
We begin this section by recalling the definition due to Yang and Yu \cite{YangYu} of \emph{lattice polarizable cubic fourfolds}:
\begin{definition}
 Let $M$ be a positive-definite lattice of rank $2\leq r(M)\leq21$.   A cubic fourfold \(X\) is called \emph{$M$-polarizable} if there exists a primitive embedding \(i:M \into A(X)\) with \(h_X^2\in i(M)\).
\end{definition}
In \cite{YangYu}, they denote the subset in $\CC$ of all $M$-polarizable cubic fourfolds by \(\mathcal{C}_M\) and prove the following theorem about these loci:\par
\begin{Thm}[{\cite[Theorem 5.1]{YangYu}}]\label{theorem:YangYu}
    Let $M$ be a positive-definite lattice of rank $r(M)\geq2$. If the following conditions hold
    \begin{enumerate}
        \item $M$ contains a \emph{distinguished element} $\o$ (that is, $(\o.\o)=3$ and the sublattice $\o^\perp\subset M$ is even),
        \item $M$ has no \emph{roots} (that is, vectors of norm 2), and
        \item $r(M)+l(M)\leq20$ (note that this is immediate for $r(M)\leq10$),
    \end{enumerate}
    then $\CC_M$ is a non-empty irreducible closed subvariety of codimension $r(M)-1$.
\end{Thm}
\begin{example}
For any natural $m\geq1$ we define the lattices:
$$\mathcal{K}_{6m+2} =\begin{pmatrix}
  & u & v \\
u & 3 & 1 \\
v & 1 & 2m+1\\
\end{pmatrix},\ \mathcal{K}_{6m} =\begin{pmatrix}
  & u & v \\
u & 3 & 0 \\
v & 0 & 2m\\
\end{pmatrix}$$
One can check that these lattices satisfy conditions (1)--(3) of \cref{theorem:YangYu} with $\o=u$, so the theorem implies that $\mathcal{C}_{\mathcal{K}_{6n}}$ and $\mathcal{C}_{\mathcal{K}_{6n+2}}$ are non-empty irreducible closed subvarieties of codimension one.  Letting $d=6n$ or $6n+2$, $\mathcal{C}_{\mathcal{K}_{d}}$ is precisely Hassett's divisor $\CC_d$ of discriminant $d$ special cubic fourfolds.

So we can see that the divisors $\mathcal{C}_d$ are a special case of $\mathcal{C}_M$, for lattices $M$ of rank $2$ and discriminant $d$ generated by the distinguished element $h^2$ and the class of an algebraic surface.
\end{example}

\subsection{Intersections of $\CC_M$}
A useful property of these subvarieties of $\CC$ is the following result:
\begin{proposition}\label{proposition:useful}
    For any cubic fourfold \(X\) and positive-definite lattice $M$ of rank $2\leq r(M)\le 21$,
\[X \in \mathcal{C}_M \iff \mathcal{C}_{A(X)}  \subset \mathcal{C}_M.\]
\end{proposition}
\begin{proof}
Suppose that $X\in \CC_M$, and let \(Y \in \mathcal{C}_{A(X)}\). Then we have a primitive isometric embedding \[\phi:A(X) \hookrightarrow A(Y)\]
such that
\(h_Y^2 \in \im(\phi)\).
As \(X \in \mathcal{C}_M\), it follows that there is another primitive isometric embedding \[\psi:M \into A(X)\] such that \(h_X^2 \in \im(\psi)\).
We know from definitions that \(h_X^2\) is a distinguished element of \(A(X)\). 
Since $\phi$ is isometric, \(\phi^{-1}(h_Y^2) \in A(X)\) has norm 3. Since $h_Y^2$ is a distinguished element in $A(Y)$, by definition we have
\[{\left(\Z h_Y^2\right)}^{\perp} \subset A(Y)\] is even. It follows that
\[{\left(\Z \phi^{-1}(h_Y^2)\right)}^{\perp} \subset A(X)\] is also even. 
Thus, we conclude that
\(\phi^{-1}(h_Y^2) \in A(X)\) is also distinguished element. By \cite[Lemma 4.8]{YangYu}, we can conclude that there exists an isometry \(\theta\) of \(A(X)\) which satisfies 
\[\theta(h_X^2)=\phi^{-1}(h_Y^2)\] 
Now, we will define \(\rho:=\phi\circ\theta\circ\psi:M \rightarrow A(Y)\). We know that
$$\rho(M)=\phi(\theta(\psi(M)))\subset\phi(\theta(A(X))=\phi(A(X))\subset A(Y),$$
where all these inclusions are primitive isometric embeddings. Thus $\rho:M\to A(Y)$ is as well.
As 
\[\rho(\psi^{-1}(h_X^2))=\phi(\theta(\psi(\psi^{-1}(h_X^2))))=\phi(\theta(h_X^2))=\phi(\phi^{-1}(h_Y^2))=h_Y^2\]
we see that \(\rho\colon M\into A(Y)\) is a primitive isometric embedding such that \(h_Y^2\in \im(\rho)\). Thus \(Y\in \mathcal{C}_M\).\par

For the converse direction, we certainly have $X\in\CC_{A(X)}$, so if $\CC_{A(X)}\subset\CC_M$ then $X\in\CC_M$.
\end{proof}

This result has an interesting consequence for the structure of intersections (finite or infinite) of subvarieties of the form $\CC_M$:
\begin{corollary}\label{corollary:intersections}
For any set of lattices \(\{M_\alpha\}_{\alpha\in A}\) the intersection of the corresponding loci in $\CC$ satisfies
\[\bigcap_{\alpha\in A} \mathcal{C}_{M_\alpha}=\bigcup_{\beta\in B}\CC_{K_\beta},\]
for a (possibly different) set of lattices \(\{K_\beta\}_{\beta\in B}\) .
\end{corollary}
\begin{proof}
    Let's compare the varieties
    \[\mathcal{I}:=\bigcap_\alpha \mathcal{C}_{M_\alpha}\]
    and
    \[\mathcal{J}:=\bigcup_{[X]\in\mathcal{I}}\mathcal{C}_{A(X)},\]
    which we claim are equal.  This suffices to prove the result.\par
    Since $[X]\in\CC_{A(X)}$, for any cubic fourfold $X$ with $[X]\in\II$, 
    $$[X]\in\CC_{A(X)}\subset\JJ,$$
    it follows that \(\mathcal{I}\subset\mathcal{J}\). \par
    For the other inclusion, let's take a cubic fourfold \([Y]\in\mathcal{J}\):
    \[[Y]\in\mathcal{J} \Rightarrow [Y]\in\bigcup_{[X]\in\mathcal{I}}\mathcal{C}_{A(X)}\Rightarrow\exists [X]\in\mathcal{I}:[Y]\in\mathcal{C}_{A(X)}.\]
    By \cref{proposition:useful} we conclude firstly that
    \[[Y] \in \mathcal{C}_{A(Y)}  \subset \mathcal{C}_{A(X)}, \]
    and secondly that $[X]\in\II$ implies $\CC_{A(X)}\subset \CC_{M_\alpha}$ for all $\alpha\in A$.
    Putting these together, it follows that $[Y]\in\II$, as required.
\end{proof}

The main example of such intersections that concerns us the intersection
$$\mathcal{Z}:=\bigcap_{d\in\HH}\ \mathcal{C}_d$$
of all Hassett divisors.
Recalling that $[X]$ being in $\mathcal{C}_d$ means that $X$ has additional geometric structure, we expect a cubic fourfold $X$ with \([X]\in\mathcal{Z}\) to have a particularly rich geometry.  
In-fact, any such $X$ is rational by the inclusions of \ZZ\ in $\CC_d$ for $d=14,26,38,42$, as it is known that for such $d$ any $X\in\CC_d$ is rational.  
Moreover, the condition of $X\in\CC_d$ for \textbf{admissible $d$} is conjectured to be equivalent to $X$ being rational, so for $[X]\in\ZZ$, we can view $X$ as rational in many ways  (conjectured to be infinitely many). These correspond to different birational maps between $X$ and $\P^4$, or equivalently, many birational automorphisms of $\P^4$. In that way 
, the rich geometry of $X\in\ZZ$ can shed light on the Cremona group of \(\P^4\).

The main result of this short note is that $\dim\ZZ=16$, answering \cite[Question 7.14]{YangYu}.
\section{The dimension of $\ZZ$}
The dimension of a variety is the maximum of the dimensions of its irreducible components. \cref{corollary:intersections} implies that any irreducible component of \SS\ and \ZZ\ is of the form $\CC_M$ for a positive-definite lattice $M$ satisfying the conditions of \cref{theorem:YangYu}. 
Now, $\CC_M\subset\CC_d$ if and only if there exist a primitive $v\in M$ such that $3\|v\|_M-\langle \o,v\rangle_M^2=\disc(\Span(\o,v))=d$. We recognize this as the quadratic form on $D(M)=M/\Z\o$ given by $\|v+\Z\o\|_{D(M)}=\disc(\Span(\o,v))$. Thus, we can conclude that for any subset $A\subset\HH$, we have 
$$\CC_M\subset\bigcap_{d\in A} \CC_d\iff A\subset \Pim(D(M)).$$
In particular, for $\ZZ=\bigcap_{d\in \HH} \CC_d$ we have: $$\CC_M\subset\ZZ\iff \HH\subset\Pim(D(M)).$$
As $\Pim(D(M))\subset \im(D(M))\subset\HH$, this inclusion must be an equality, so we have the following result:
\begin{proposition}\label{Prop:Codim=rank}
    $$\codim(\ZZ) = \min\left\{r(M)-1|\CC_M\subset\ZZ\right\} =  \min\left\{r(D(M))| \HH = \Pim(D(M))\right\}=r_{prim}(\HH).$$
\end{proposition}

In \cref{corollary:RankHFour} we computed that $r_{prim}(\HH)=4$, so we can conclude that:
\begin{Thm}\label{Thm:MainTheorem}
    $\dim(\mathcal{Z})=16$.
\end{Thm}
\section{Appendix: Some auxiliary arithmetic results}
In determining the primitive supporting rank of $\HH$ we used the following arithmetic results, the first of which is interesting in its own right.
\begin{proposition}
For any natural $n\not\equiv7\cdot4^k \pmod {8\cdot4^k}$ there exists a non-trivial solution to
\begin{equation}\label{eq:abcsquares}
    a^2+b^2+c^2=n
\end{equation}
 if and only if $n\not=4^k$.  Here a solution $(a,b,c)$ is called trivial if $(ab,bc,ca)=(0,0,0)$, that is, two out of three of the coordinates vanish.
\end{proposition}
\begin{proof}
By Legendre's Theorem on sums of three squares, any $n$ as in the hypothesis admits some representation as a sum of three squares, so the question at hand is whether a non-trivial representation can be chosen.  If $n$ not a perfect square, then any representation is automatically non-trivial, so may assume that \(n\) is a perfect square. 

Moreover, multiplying a non-trivial solution $(a,b,c)$ to \eqref{eq:abcsquares} for \(n=k^2\) by \(m\) gives the non-trivial $(ma,mb,mc)$ to \eqref{eq:abcsquares} for $m^2k^2$.  Thus, it suffices to treat the case \(n=p^2\) for a prime $p$, which may be taken to be odd since $n\ne 4^k$.
If \(p\equiv1\pmod 4\), then by Fermat's theorem on sums of two squares it follows that we can write \(p=x^2+y^2\).  Therefore 
$$p^2=(x^2+y^2)^2=(x^2-y^2)^2+(2xy)^2=(x^2-y^2)^2+(2xy)^2+0^2.$$
This is in-fact a non-trivial representation because $2xy\ne 0$, since $p$ is prime and thus not a perfect square, and $x^2-y^2\ne0$, since then $p=2$.

It remains to treat the case \(p\equiv3\pmod 4\).
Consider the equation 
\begin{equation}\label{eq:abcsquares4p^2}
    a^2+b^2+c^2=4p^2
\end{equation}
and observe by reducing $\mod 4$ that any solution $(a,b,c)$ to \eqref{eq:abcsquares4p^2} must lie in $(2\N)^3$.  Thus if we show that there exists a non-trivial solution $(a,b,c)$ to \eqref{eq:abcsquares4p^2}, then we can divide by $2$ to obtain the non-trivial solution 
\[\left(\frac{a}{2}\right)^2+\left(\frac{b}{2}\right)^2+\left(\frac{c}{2}\right)^2=p^2,\]
as required.

As $p$ is an odd prime, by Lagrange's three-square theorem \(2p\) can be written as a sum of 3 squares, say
\[2p=s^2+t^2+u^2.\]
Then 
\[4p^2=(s^2+t^2+u^2)^2=(t^2+u^2-s^2)^2+(2st)^2+(2su)^2\]
is a representation as a sum of three squares, and we claim it is a nontrivial representation.  Indeed, if we assume to the contrary then we conclude either that $p$ or $2p$ is a perfect square, which is impossible, or that $2p$  is a sum of two squares, contrary to Fermat's two-square theorem and the fact that $p\equiv 3\pmod 4$.  This finishes the proof of one direction.

For the converse, suppose that $n=4^k$ and consider a representation as a sum of three squares
$$4^k=a^2+b^2+c^2.$$
As above, we must have that $(a,b,c)\in(2\N)^3$, so dividing the solution by $2$ gives
$$4^{k-1}=\left(\frac{a}{2}\right)^2+\left(\frac{b}{2}\right)^2+\left(\frac{c}{2}\right)^2,$$
and continuing in this fashion, we obtain 
$$1=\left(\frac{a}{2^k}\right)^2+\left(\frac{b}{2^k}\right)^2+\left(\frac{c}{2^k}\right)^2.$$
From this it is clear that the solution must be trivial.
\end{proof}
\begin{lemma}\label{lemma:Lagrange5}
$$x^2+y^2=5n\text{ has an integral solution} \iff x^2+y^2=n\text{ has an integral solution}.$$
\end{lemma}
\begin{proof}
By Fermat's two-square theorem  \(x^2+y^2=n\) has an integral solution if and only if any prime divisor \(p\) of \(n\) with \(2\nmid r_p(n)\) satisfy \(p\equiv1,2\pmod 4\), where $r_p(n)$ is the biggest integer \(r\) such \(p^r|n\).  It is then immediate that \(5n\) satisfies this condition if and only if \(n\) satisfies this condition.
\end{proof}
\begin{proposition}\label{proposition:Lowering the 5-ht}
If $x^2+y^2=5n\not=0$ is solvable over $\Z$, then there exists a solution which satisfies $x^2\equiv-y^2\equiv1\pmod 5$.
\end{proposition}
\begin{proof}
Write $n$ as $n=5^km$ for $k\geq0,\ m\not\equiv0\pmod 5$. We will prove the proposition by induction on $k$.\par
\textbf{Base}: For $k=0$, we have integers $x,y\in\Z$ such that $x^2+y^2=5m$.  
It follows that $x^2\equiv-y^2\pmod 5$. 
Up to switching $x$ and $y$, we can assume that $x^2\equiv-y^2\equiv1\pmod 5$ or $x^2\equiv-y^2\equiv0\pmod 5$. 
In the second case, it would follow that $x\equiv y\equiv0\pmod 5$, so in fact 
$$x^2\equiv y^2\equiv0\pmod {25}\Longrightarrow 5m=x^2+y^2\equiv 0\pmod {25}.$$
But then $m\equiv0\pmod 5$, contradicting $m\not\equiv0\pmod 5$.\par
\textbf{Induction Step:} Let's assume that the proposition holds for $k-1$, and take $n=5^km$. \cref{lemma:Lagrange5} implies that 
$$x^2+y^2=5^{k+1}m=5n\text{ has an integral solution} \iff x^2+y^2=5^km=n\text{ has an integral solution}.$$
So, if $x^2+y^2=5n$ is solvable, then $x^2+y^2=5^km$ is solvable, and by the induction hypothesis we can choose integers $(a,b)$ such that $(x,y)=(5a+1,5b+2)$ satisfy
$(5a+1)^2+(5b+2)^2=5^km$. Now we compute 
$$\left(10a+5b+4\right)^2+\left(10b-5a+3\right)^2=5((5a+1)^2+(5b+2)^2)=5^{k+1}m=5n$$
to demonstrate the existence of a solution satisfying the desired property.
\end{proof}
\begin{proposition}\label{proposition:Mod3Solutions}
    If \(x^2+5y^2=2n\) for \(n\not=0\) is solvable over $\Z$, then there exists a solution which satisfies \(x,y\not\equiv0\pmod 3\). Moreover, for any pair of solutions \((x_1,y_1),(x_2,y_2)\), we have
    \[x_1\equiv y_1\equiv x_2\equiv y_2\pmod 2,\qquad\ x_1\equiv\pm x_2\pmod 5.\]
\end{proposition}
\begin{proof}
For the first statement, we can use the same method of proof as in \cref{proposition:Lowering the 5-ht}, with the help of the equation
    \[x^2+5y^2=\frac{1}{9}\left(2^2+5\cdot1^2\right)\left(x^2+5y^2\right)=\left(\frac{2x\pm5y}{3}\right)^2+5\left(\frac{x\mp2y}{3}\right)^2\]
    which helps us to reduce the 3-height of $x$ amd $y$. \par
    For the second statement, reducing the equation  \(x^2+5y^2=2n\) modulo $2$ and $4$ implies that \(x_i\equiv y_i\equiv n\pmod 2\). Moreover, \(x_1^2\equiv x_2^2\equiv2n\pmod 5\) which implies that $x_1\equiv\pm x_2\pmod 5$, as we want.
    
\end{proof}

\bibliographystyle{plain}
\bibliography{main}
\end{document}